# A novel method in solving seepage problems implementation in Abaqus based on the polygonal scaled boundary finite element method

Yang Yang [1,2,*], Zongliang Zhang [1,*], Yelin Feng [1]

1. PowerChina Kunming Engineering Corporation Limited, Kunming 650051, China ;
2. Department of Hydraulic Engineering, Tsinghua University, Beijing 100084, China;
* Correspondence: Yang Yang, yangyhhu@foxmail.com;

**Abstract**

The scaled boundary finite element method (SBFEM) is a semi-analytical computational scheme, which is based on the characteristics of the finite element method (FEM) and combines the advantages of the boundary element method (BEM). This paper integrates the scaled boundary finite element method (SBFEM) and the polygonal mesh technique into a new approach to solving the steady-state and transient seepage problems. The proposed method is implemented in Abaqus using a user-defined element (UEL). The detailed implementations of the procedure, defining the UEL element, updating the RHS and AMATRX, and solving the stiffness/mass matrix by the eigenvalue decomposition are presented. Several benchmark problems from seepage are solved to validate the proposed implementation. Results show that the polygonal element of PS-SBFEM has a higher accuracy rate than the standard FEM element in the same element size. For the transient problems, the results between PS-SBFEM and the FEM are in excellent agreement. Furthermore, the PS-SBFEM with quadtree meshes shows a good effect for solving complex geometric in the seepage problem. Hence, the proposed method is robust accurate for solving the steady-state and transient seepage problems. The developed UEL source code and the associated input files can be downloaded from https://github.com/yangyLab/PS-SBFEM.

## 1. Introduction

Seepage analysis is an essential topic in civil engineering. The changes in soil pore water pressure may significantly affect the stability of structures, such as slope engineering [1], tunnel engineering [2], and earth-rock dam engineering [3]. The finite element method (FEM) is one of the dominant methods for seepage problems [4][5][6]. However, this method is cumbersome in dealing with the singularities. Furthermore, the whole computational domain must be discretized in the FEM [7].

Regarding the shortcomings, an alternative approaches have been proposed. The scaled boundary finite element method (SBFEM) was developed in the 1990s. The scaled boundary finite element method is a semi-analytical method that attempts to combine the advantages and characteristics of FEM and the boundary element method (BEM) into one new approach. In the SBFEM, only the boundaries of the domain are discretized in the circumferential direction. Then, in the radial direction, the partial differential equation (PDE) is transformed to an ordinary differential equation (ODE) which can be solved analytically [8]. The SBFEM has been applied to many engineering problems, such as wave propagation [9,10], heat conduction [11,12], fracture [13–16], acoustic [17], seepage [18,19], elastoplastic [20], fluid [21,22].

The polygonal scaled boundary finite element method (PSBFEM) is a novel method integrating the standard SBFEM and the polygonal mesh technique [20,23]. Compared with the standard finite element method, the polygon element with more than four edges involves more nodes in the domain and usually shows higher accuracy [20]. Polygon can discretize complex geometry flexibly. Furthermore, polygons have high geometric isotropy and eliminate the mesh dependence caused by the discretization of fixed meshes with standard triangles or quadrangles [24]. These advantages further motivate polygonal finite elements as an alternative to standard finite elements using triangles or quadrangles.

Recently, an alternative mesh technique has been widely used in geometric discretization. The quadtree algorithm is fast, efficient, and capable of achieving rapid and smooth transitions of element sizes between mesh refinement regions [25]. Due to hanging nodes between two adjacent elements of different sizes, it is problematic that quadtree meshes are directly used to simulate within the finite element method. However, the SBFEM only discretizes in the boundary of geometry. Hence, each element in a quadtree mesh is treated as a generic polygon regardless of hanging nodes. The enables the structure of the quadtree to be exploited for efficient computations. The ability to assume any number of sides also enables the SBFEM to discretize curved boundaries better.

Until now, few studies of SBFEM seepage analysis have been reported. Li and Tu [24] used the SBFEM to solve the steady-state seepage problems with multi-material regions. Bazyar and Talebi [26] simulated the transient seepage problems in zoned anisotropic soils. Prempramote [27] developed a high-frequency open boundary for the transient seepage analyses of semi-infinite layers with a constant depth. Liu et al. [18] presented an iso-geometric scaled boundary finite element method (IGA-SBFEM) using the non-uniform rational B-splines (NURBS) for the numerical solution of seepage problems in the unbounded domain. These studies demonstrate that SBFEM has excellent accuracy, efficiency, and convergence rate. However, the PSBFEM has not been used to solving the seepage problem.

Although SBFEM has matured, these studies only exist as independent codes, and SBFEM has not become commercial software. Therefore, engineers do not find it easy to use this method to solve engineering problems. The commercial software Abaqus has powerful linear or nonlinear, static, or dynamic analysis capabilities [28]. Abaqus/Standard analysis also provides a user-defined element (UEL) to define an element with an available interface to the code. Several works so far have focused on the implementation of SBFEM in the Abaqus. Yang el al. [29] presents available UEL subroutines of a steady-state and transient heat conduction analysis using the PSBFEM. Ya et al. [30] implemented an open-source polyhedral SBFEM element for three-dimensional and nonlinear problems through the Abaqus UEL. Yang et al. [31] implemented the SBFEM in Abaqus in linear elastic stress analyses. However, there are no available subroutines of SBFEM to solve the seepage problem in the Abaqus nowadays.

The primary goal of this paper is to implement a novel semi-analytical approach by integrating the scaled boundary finite element method and polygonal mesh technique to solving the seepage problem. This paper is divided into six sections: Section 2 describes the polygon seepage scaled boundary finite element method (PS-SBFEM) concept. Section 3 outlines the solution procedures. Section 4 describes the implementation of PS-SBFEM in the seepage problems by the Abaqus UEL subroutine. Then, several benchmark examples are solved in Section 5. Finally, the main concluding remarks of this paper are presented in Section 6.

## 2. The PS-SBFEM for the transient seepage problem

The governing equations related to two-dimensional transient seepage flow can be written as [19]

$$\nabla \cdot (k \nabla h) + p - S_s \dot{h} = 0 \quad (1)$$

where $S_s$ is a specific storage coefficient, $p$ is the source per unit volume, and $h$ is the total head, $\dot{h}$ is the derivative of total head with respect to time, $k$ is the permeability matrix, $\nabla$ is the gradient operator, $\nabla = \frac{\partial}{\partial x}\vec{i} + \frac{\partial}{\partial y}\vec{j}$. By applying the Fourier transform on the governing equation [26,32], it is transformed into the frequency domain as

$$\nabla \cdot (k \nabla \tilde{h}) + \tilde{p} - i\omega S_s \tilde{h} = 0 \quad (2)$$

where $\tilde{p}$, and $\tilde{h}$ are the Fourier transform of $p$, and $h$, respectively. $\omega$ is the frequency.

As illustrated in Figure 1, the SBFEM presents a local coordinate system $(\xi,\eta)$. The coordinates of any point $(\tilde{x}, \tilde{y})$ along the radial line and inside the domain can be written as [8]

$$\tilde{x} = \xi[N(\eta)]\{x\} \tag{3}$$

$$\tilde{y} = \xi[N(\eta)]\{y\} \tag{4}$$

where $[N(\eta)]$ is the shape function matrix.

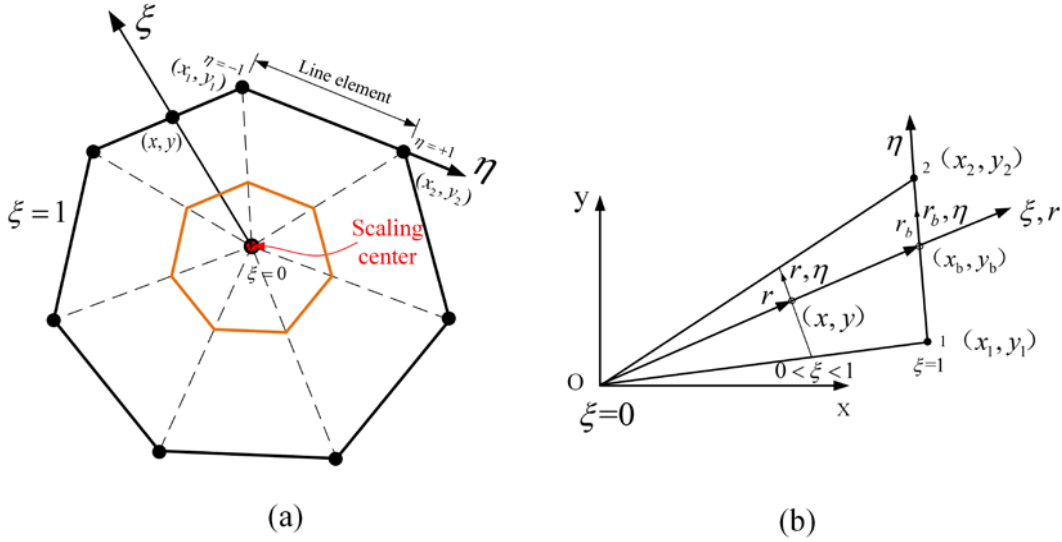

Figure 1. The coordinate system of SBFEM; (a) the $S$-element; (b) the local coordinate system of SBFEM.

The differential operator in the Cartesian coordinate system can be transformed to the scaled boundary coordinate system as follows [8]:

$$\nabla = [b_1]\frac{\partial}{\partial \xi} + \frac{1}{\xi}[b_2]\frac{\partial}{\partial \eta} \tag{5}$$

with

$$[b_1] = \frac{1}{|J_b|}\begin{bmatrix} y_{b,\eta} & 0 \\ 0 & -x_{b,\eta} \end{bmatrix} \tag{6}$$

$$[b_2] = \frac{1}{|J_b|}\begin{bmatrix} -y_b & 0 \\ 0 & x_b \end{bmatrix} \tag{7}$$

where the Jacobian matrix at the boundary can be written as

$$[J_b] = \begin{bmatrix} x_b & y_b \\ x_{b,\eta} & y_{b,\eta} \end{bmatrix} = x_b y_{b,\eta} - y_b x_{b,\eta} \tag{8}$$

The head function at any point can be expressed as

$$\{\tilde{h}(\xi,\eta)\} = [N_u(\eta)]\{\tilde{h}(\xi)\} \tag{9}$$

where $\{\tilde{h}(\xi)\}$ is the nodal head vector. $[N_u(\eta)]$ is the shape function matrix.

By using Equation (5) and (9), the flux $\tilde{Q}(\xi,\eta)$ can be written as

$$\tilde{Q}(\xi,\eta) = -[k]\left([B_1(\eta)]\{\tilde{h}(\xi)\}_{,\xi} + \frac{1}{\xi}[B_2(\eta)]\{\tilde{h}(\xi)\}\right) \tag{10}$$

with

$$[B_1(\eta)] = \{b_1(\eta)\}[N(\eta)] \tag{11}$$

$$[B_2(\eta)] = \{b_2(\eta)\}[N(\eta)] \tag{12}$$

Applying the weighted residual method and Green's theorem and introducing the boundary conditions, the following equations can be written as follows [8,21]

$$[E_0]\xi^2\{\tilde{h}(\xi)\}_{,\xi\xi} + \left([E_0]-[E_1]+[E_1]^T\right)\xi\{\tilde{h}(\xi)\}_{,\xi} - \left([E_2]-i\omega[M_0]\xi^2\right)\{\tilde{h}(\xi)\} = \xi\{F(\xi)\} \tag{13}$$

where,

$$[E_0] = \int_S [B_1(\eta)]^T k[B_1(\eta)]|J_b|d\eta \tag{14}$$

$$[E_1] = \int_S [B_2(\eta)]^T k[B_1(\eta)]|J_b|d\eta \tag{15}$$

$$[E_2] = \int_S [B_2(\eta)]^T k[B_2(\eta)]|J_b|d\eta \tag{16}$$

$$[M_0] = \int_S [N(\eta)]^T S_s[N(\eta)]|J_b|d\eta \tag{17}$$

## 3. Solution procedure of the PS-SBFEM equation

3.1 Steady-state solution

Equation (13) is an second-order Euler-Cauchy equation. In this paper, considering the side is adiabatic or the domain is closed, the $\{F(\xi)\}$ on the right side of Equation (13) satisfies $\{F(\xi)\}=0$. Setting $\omega=0$ in Equation (13), the SBFEM equations for the steady-state seepage field can be written as

$$[E_0]\xi^2\{\tilde{h}(\xi)\}_{,\xi\xi} + \left([E_0]+[E_1]^T-[E_1]\right)\xi\{\tilde{h}(\xi)\}_{,\xi} - [E_2]\{\tilde{h}(\xi)\} = 0 \tag{18}$$

By introducing the variable $\{X(\xi)\}$ that consists of nodal water head functions $\tilde{h}(\xi)$ and the flux functions $\tilde{Q}(\xi)$,

$$\{X(\xi)\} = \begin{Bmatrix} \{\tilde{h}(\xi)\} \\ \{\tilde{Q}(\xi)\} \end{Bmatrix} \tag{19}$$

The equation can be transformed into the first-order ordinary differential equation:

$$\xi\{X(\xi)\},\xi - [Z_p]\{X(\xi)\} = 0 \tag{20}$$

where the coefficient matrix $[Z_p]$ is a Hamiltonian matrix. The solution for a bounded domain is obtained using the positive eigenvalues of $[Z_p]$. Hence, $[Z_p]$ can be expressed as

$$[Z_p] = \begin{bmatrix} -[E_0]^{-1}[E_1]^T & [E_0]^{-1} \\ [E_2]-[E_1][E_0]^{-1}[E_1]^T & [E_1][E_0]^{-1} \end{bmatrix} \tag{21}$$

The solution of Equation (21) can be obtained by computing the eigenvalue and eigenvector of the matrix $[Z_p]$ yields

$$[Z_p]\begin{bmatrix} [\psi_{11}] & [\psi_{12}] \\ [\psi_{21}] & [\psi_{22}] \end{bmatrix} = \begin{bmatrix} [\psi_{11}] & [\psi_{12}] \\ [\psi_{21}] & [\psi_{22}] \end{bmatrix}\begin{bmatrix} [\lambda_n] & \\ & [\lambda_p] \end{bmatrix} \tag{22}$$

where the real components of eigenvalues $\lambda_n$ and $\lambda_p$ are negative and positive, respectively.

The general solution of equation (21) can be obtained as follows:

$$\{\tilde{h}(\xi)\} = [\psi_{11}]\xi^{-[\lambda_1]}\{c_1\} + [\psi_{12}]\xi^{-[\lambda_2]}\{c_2\} \tag{23}$$

$$\{\tilde{Q}(\xi)\} = [\psi_{21}]\xi^{-[\lambda_1]}\{c_1\} + [\psi_{22}]\xi^{-[\lambda_2]}\{c_2\} \tag{24}$$

To obtain a finite solution at the scaling center ($\xi = 0$), $\{c_2\}$ must be equal to zero, the solution in bounded domain can be written as:

$$\{\tilde{h}(\xi)\} = [\psi_{11}]\xi^{-[\lambda_1]}\{c_1\} \tag{25}$$

$$\{\tilde{Q}(\xi)\} = [\psi_{21}]\xi^{-[\lambda_1]}\{c_1\} \tag{26}$$

The relationship between $\{\tilde{h}(\xi)\}$ and $\{\tilde{Q}(\xi)\}$ is expressed as

$$[K^{st}]\{\tilde{h}(\xi)\} = \{\tilde{Q}(\xi)\} \tag{27}$$

where the steady-state stiffness matrix of the S-element can be expressed as

$$[K^{st}] = [\psi_{21}][\psi_{11}]^{-1} \tag{28}$$

## 3.2 Mass matrix and transient solution

To determine the mass matrix $[M]$ of SBFEM, introducing the dynamic-stiffness matrix $[K(\xi,\omega)]$ at $\xi$:

$$[K(\xi,\omega)]\{\tilde{h}(\xi)\} = \{\tilde{Q}(\xi)\} \tag{29}$$

For the bounded domain, the dynamic-stiffness matrix $[K(\omega)]$ on the boundary $\xi=1$ formulated in the frequency domain is written as

$$\left([K(\omega)]-[E_1]\right)[E_0]^{-1}\left([K(\omega)]-[E_1]^T\right)-[E_2]+2\omega[K(\omega)]_{,\omega}-i\omega[M_0]=0 \tag{30}$$

To obtain the mass matrix $[M]$ of the bound domain, the low-frequency case in the SBFEM is address, which the dynamic-stiffness $[K(\omega)]$ of a bounded domain is assume as

$$[K(\omega)] = \left[K^{st}\right] + i\omega[M] \tag{31}$$

The steady-state stiffness matrix $[K^{st}]$ is computed using Equation (28). Substituting Equation (31) into Equation (30) leads to a constant term independent of $i\omega$, a term in $i\omega$ and higher order term in $i\omega$, which are neglected. Besides, the constant term vanishes. The coefficient matrix of $i\omega$ can be expressed as

$$\left(\left(-\left[K^{st}\right]+[E_1]\right)[E_0]^{-1}-[I]\right)[M]+[M]\left([E_0]^{-1}\left(-\left[K^{st}\right]+[E_1]^T\right)-[I]\right)+[M_0]=0 \tag{32}$$

This is a linear equation to solver the mass matrix $[M]$. By using the eigenvalues and eigenvectors of matrix $[Z_P]$, the Equation (32) can be written as

$$\left([I]+[\lambda_b]\right)[m]+[m]\left([I]+[\lambda_b]\right) = \left[\Phi_b^{\tilde{T}}\right]^T[M_0]\left[\Phi_b^{\tilde{T}}\right] \tag{33}$$

where

$$[m] = \left[\Phi_b^{\tilde{T}}\right]^T[M]\left[\Phi_b^{\tilde{T}}\right] \tag{34}$$

After solving matrix $[m]$ in the Equation (34), the mass matrix $[M]$ is obtained by

$$[M] = \left(\left[\Phi_b^{\tilde{T}}\right]^{-1}\right)^T[m]\left[\Phi_b^{\tilde{T}}\right]^{-1} \tag{35}$$

3.3 Time discretization

Substitute Equation (31) into Equation (29) followed by performing inverse Fourier transform, the nodal water head relationship of a bounded domain is expressed as a standard time domain equation using steady-state stiffness and mass matrices as

$$[K^{st}]\{h(t)\}+[M]\{\dot{h}(t)\}=\{Q(t)\} \tag{36}$$

where nodal water head $h(t)$ is the continuous derivative of time, In general, it is not easy to solve the function solution in the time domain. In this paper, the backward difference method [34,35] is adopted to solve Equation (36). The time domain is divided into several time units, and the solution of the time node is obtained step by step from the initial conditions, and the nodal water head at any time is obtained by interpolation.

At time $[t, t+\Delta t]$, the water head change rate $\{\dot{h}(t)\}$ can be expressed as

$$\{\dot{h}(t)\} = \frac{[\Delta h]}{\Delta t} = \frac{\{h(t)\}^{t+\Delta t} - \{h(t)\}^{t}}{\Delta t} \tag{37}$$

Substitute Equation (40) into Equation (39), the equation at time step $t+\Delta t$ can be obtain as follows

$$([K^{st}]^{t+\Delta t} + \frac{[M]^{t+\Delta t}}{\Delta t})\{h(t)\}^{t+\Delta t} = \{Q(t)\}^{t+\Delta t} + \frac{[M]^{t+\Delta t}}{\Delta t}\{h(t)\}^{t} \tag{38}$$

**4. The procedure implementation**

4.1 Implementation of PS-SBFEM in the ABAQUS

The most critical work of UEL is to update the contribution of the element to the internal force vector RHS and the stiffness matrix AMATRX in the user subroutine interface provided by Abaqus [35]. For the steady-state seepage analysis, the AMATRX and RHS is defined as following:

$$\text{AMATRX} = [K^{st}] \tag{39}$$

$$\text{RHS} = -[K^{st}]\{U\} \tag{40}$$

where $\{U\}$ is the water head vector.

For the transient seepage, the AMATRX and RHS is defined as following:

$$\text{AMATRX} = \left[K^{st}\right]^{t+\Delta t} + \frac{[M]^{t+\Delta t}}{\Delta t} \tag{41}$$

$$\text{RHS} = -\left[K^{st}\right]^{t+\Delta t} \{U\}^{t+\Delta t} - \frac{[M]^{t+\Delta t}}{\Delta t}(\{U\}^{t+\Delta t} - \{U\}^{t}) \tag{42}$$

Figure 2 illustrates the flowchart of UEL subroutine for PS-SBFEM. Based on the input file's connectivity information, the UEL computes the scaling centers and transforms the global coordinate into the local coordinate. Equations (14)~(17) computes the coefficient matrices $[E_0]$, $[E_1]$, $[E_2]$, and $[M_0]$, which are used to construct the Hamilton matrix $[Z_p]$ using the Equation (21). The two eigenvector matrices ($\left[\Phi_q^{(n)}\right], \left[\Phi_u^{(n)}\right]$) are constructed using the eigenvalue decomposition. Finally, we can obtain the stiffness matrix $[K]$ and mass matrix $[M]$ of the PS-SBFEM elements.

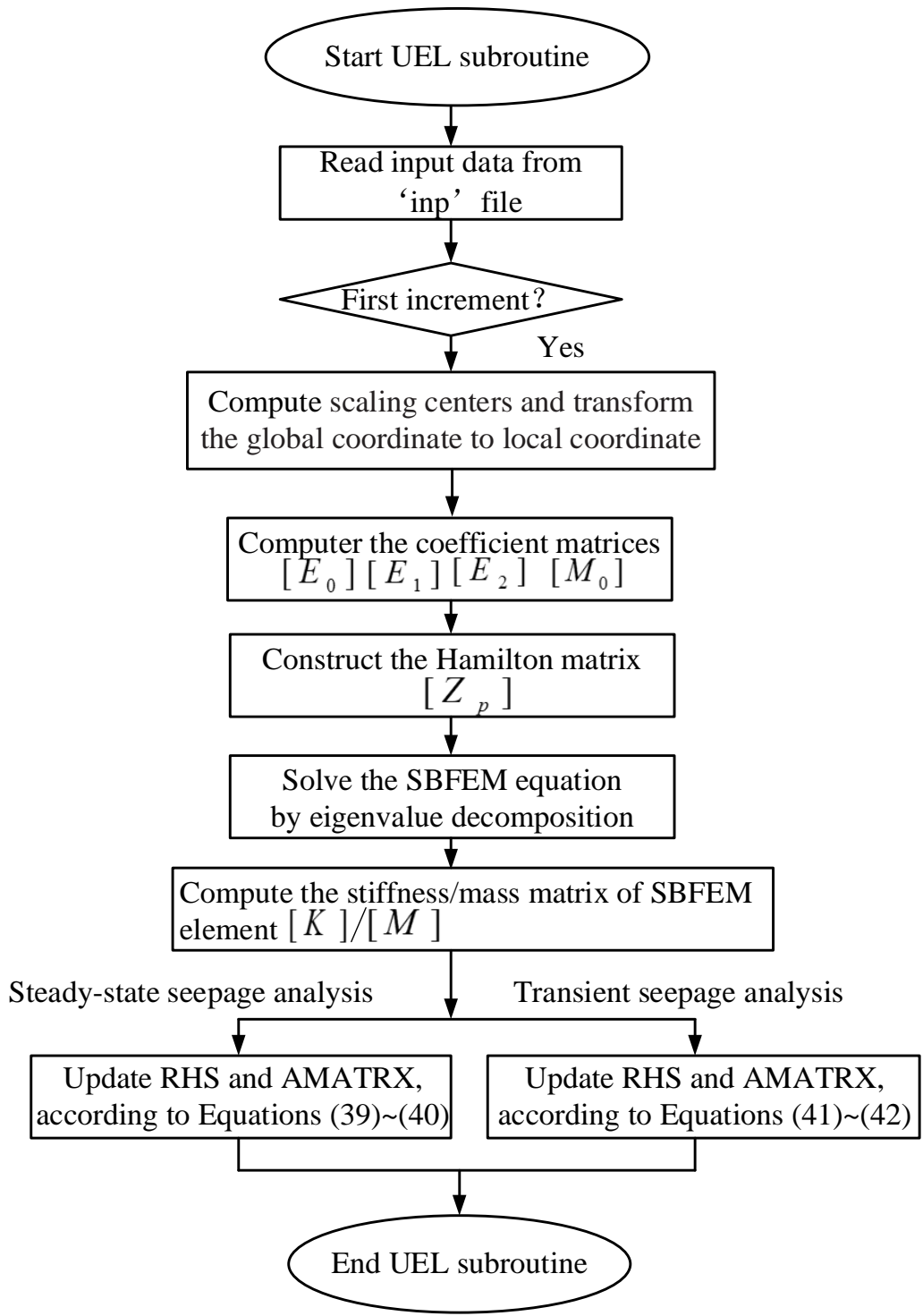

Figure 2. Flowchart of UEL subroutine for PS-SBFEM.

4.2 Polygonal mesh generation

In this work, we developed the automatic generation program of the polygonal mesh based on the polygon discretization [23,24] and the quadtree decomposition [25] . This program is developed by the Python script and the visualization in the Paraview [36]. It can be observed from Figure 3 that the polygonal mesh and quadtree mesh have an excellent ability to process complex geometric boundaries and high automation and strong universality, which can provide powerful pre-processing for engineering problems.

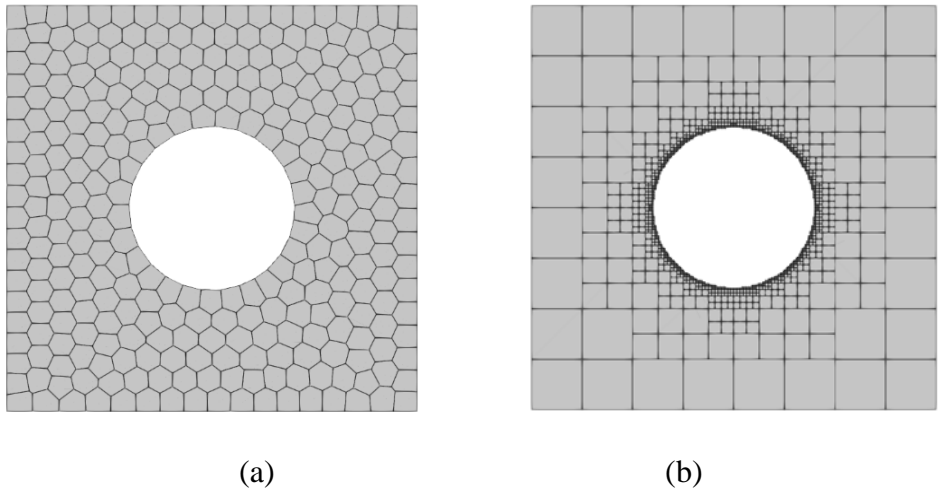

(a)  (b)

Figure 3. Polygon mesh automatically generated by Python script; (a) the polygonal mesh; (b) the quadtree mesh;

4.3 Defining the element of PS-SBFEM

The input file of Abaqus usually contains a complete description of the numerical model, such as nodes, elements, the degree of freedom, and materials. This information needs to be defined in the "inp" file by the user. We show a simple polygon mesh of PS-SBFEM to demonstrate the definition of elements in the UEL, as shown in Figure 4. The mesh consists of three element types: triangular element (U3), quadrilateral element (U4), and pentagon element (U5). As shown in Listing 1, the pentagon element (U5) is defined as follows: **1 ~ 6** is the line number. Lines **1 ~ 6** are used to define the pentagonal element (U5); Line **1** assigns the element type, the number of nodes, the number of element properties, and the number of freedom degrees per node; Line **2** sets the active degrees of freedom for the pore pressure. Lines **3~4** define the element sets' E5'; Lines **5~6** sets the permeability coefficient of 'E5'.

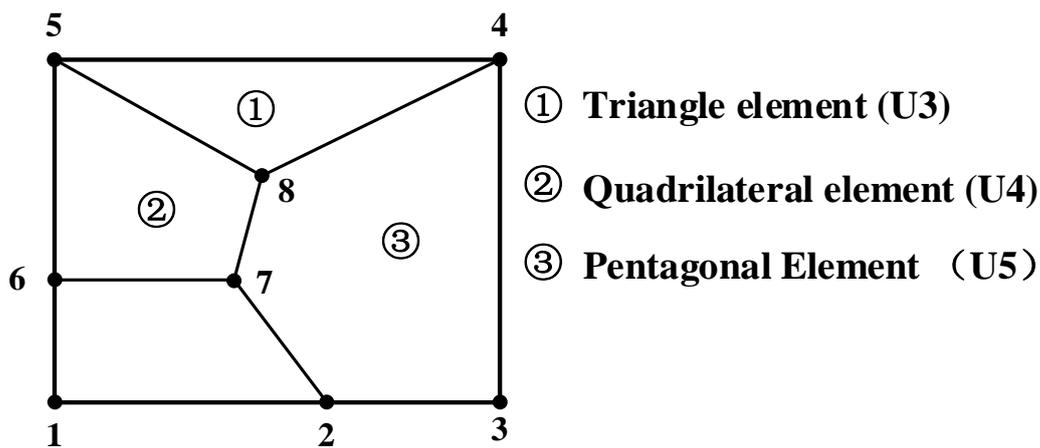

① Triangle element (U3)
② Quadrilateral element (U4)
③ Pentagonal Element （U5）

Figure 4. Schematic diagram of the polygon meshes

Listing. 1 The input file of polygon element in the Abaqus (c.f. Figure 4)

```
1 *USER ELEMENT, NODES=5, TYPE=U5, PROPERTIES=2, COORDINATES=2
2 8
3 *ELEMENT, TYPE=U5, ELSET=E5
4 3,2,3,4,8,7
```

```
5 *UEL PROPERTY, ELEST=E5
6 0.003,0.003
```

## 5. Numeric example

In this section, we carried out four numerical examples to demonstrate the convergence and accuracy of the PS-SBFEM. In addition, the results of the PS-SBFEM also are compared with the standard FEM. The FEM analysis uses the commercial finite element software Abaqus. For validation, the relative errors in the water head are investigated as follows:

$$e_{L^2} = \|\mathbf{h}\text{-}\mathbf{h}^n\|_{L^2(\Omega)} = \frac{\sqrt{\int_\Omega (\mathbf{h}-\mathbf{h}^n)^T (\mathbf{h}-\mathbf{h}^n) d\Omega}}{\sqrt{\int_\Omega \mathbf{h}^T \mathbf{h} d\Omega}} \tag{43}$$

where $\mathbf{h}$ is the numerical solution, and $\mathbf{h}^n$ is the analytical or reference solution.

5.1 Steady-state seepage problem in the concrete dam

In the first example, we consider a standard concrete dam foundation steady-state seepage problem. The geometry model and boundary conditions are shown in Figure 5 (a). The boundaries of BC, AE, EF, and DF are defined as impermeable boundaries. The boundaries of AB and CD are the water head boundary. To verify the accuracy of the proposed method, three monitor points **1** (100,80), **2** (120,80), and **3** (140,80) are chosen, as shown in Figure 5 (a). The permeability coefficient of the dam foundation is $k_x = k_y = 1\times10^{-5}$ cm/s. The FEM analysis uses the Abaqus CPE4P element in this work, as shown in Figure 5 (b). The PS-SBFEM uses the polygonal element, as shown in Figure 5 (c).

The water head of monitor points is presented in Table. 1. The relative error of FEM [37] is 4.08%; the relative error of Abaqus CPE4P and PS-SBFEM is 1.38% and 0.90%, respectively. Hence, the accuracy of PS-SBFEM is greater than FEM at the same element size. Figure 6 shows that the results are virtually the same for the FEM and PS-SBFEM. In addition, a convergence study is performed by mesh *h*-refinement. The meshes are refined successively following the sequence 20m, 10m, 5m, 2.5m, and 1.25m. Figure 7 illustrates that the water head of monitor points at different freedom degrees. It is noted that the convergence rate of PS-SBFEM and Abaqus CPE4P is the same as the mesh refinement. Moreover, the accuracy of PS-SBFEM is higher than that of Abaqus CPE4P in the same freedom degrees.

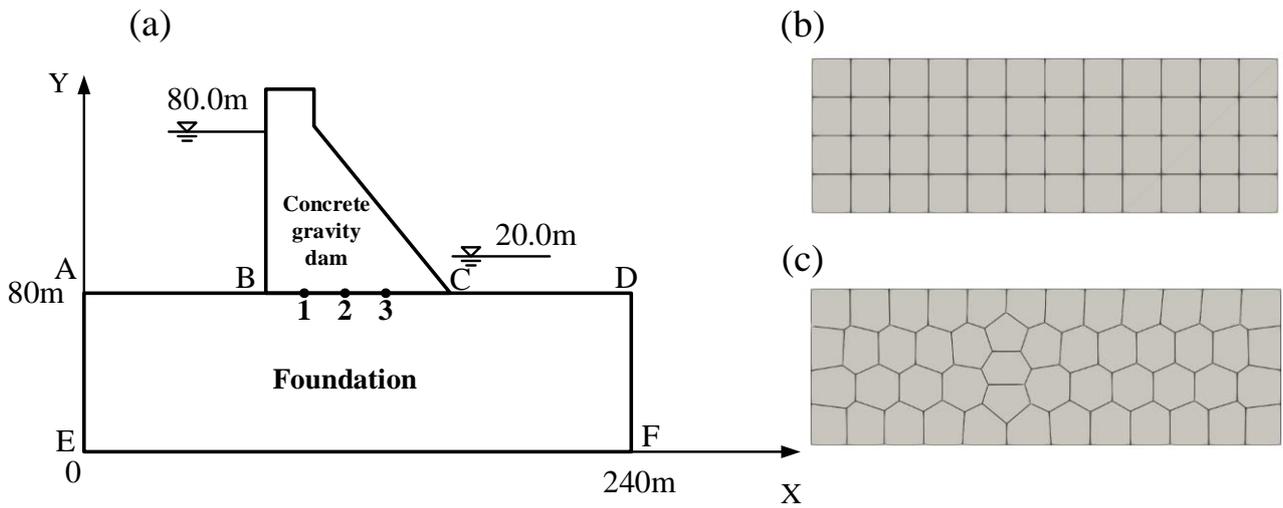

Figure 5. The steady-state seepage problem in the concrete dam; (a) geometric model and boundary conditions; (b) Abaqus CPE4P element; (c) PS-SBFEM polygonal element.

Table. 1 Comparison of water head by different methods (element size 20m)

| Method | Monitor point | | | Relative error $e_{L^2}$ (%) |
|---|---|---|---|---|
| | 1 | 2 | 3 | |
| Analytical solution (m) | 60 | 50 | 40 | - |
| FEM (m) [37] | 62.45 | 49.90 | 37.39 | 4.08 |
| Abaqus CPE4P (m) | 60.86 | 50 | 39.14 | 1.38 |
| PS-SBFEM (m) | 60.40 | 50.06 | 39.32 | 0.90 |

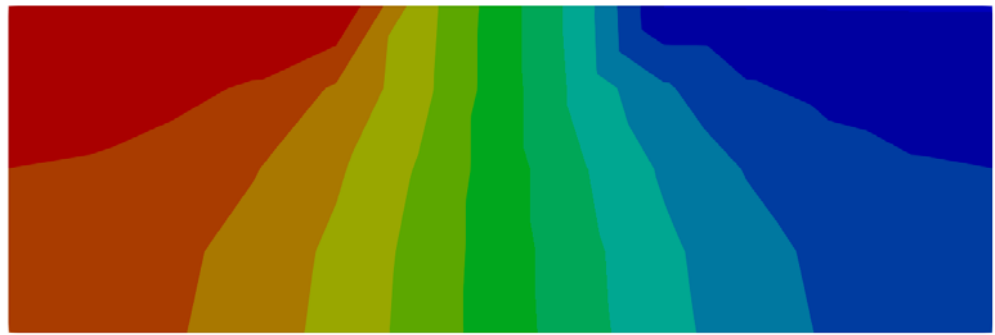
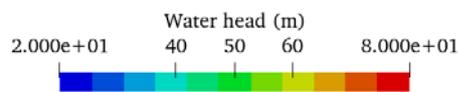

(a)

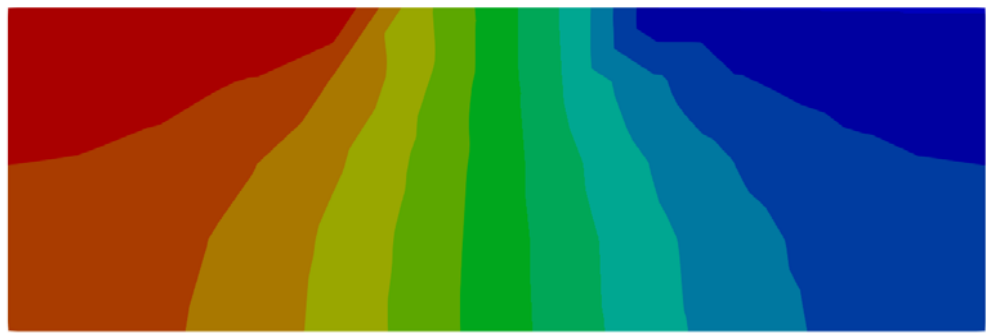
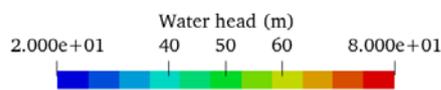

(b)

Figure 6. The water head distribution of dam foundation; (a) Abaqus CPE4P element; (b) PS-SBFEM element

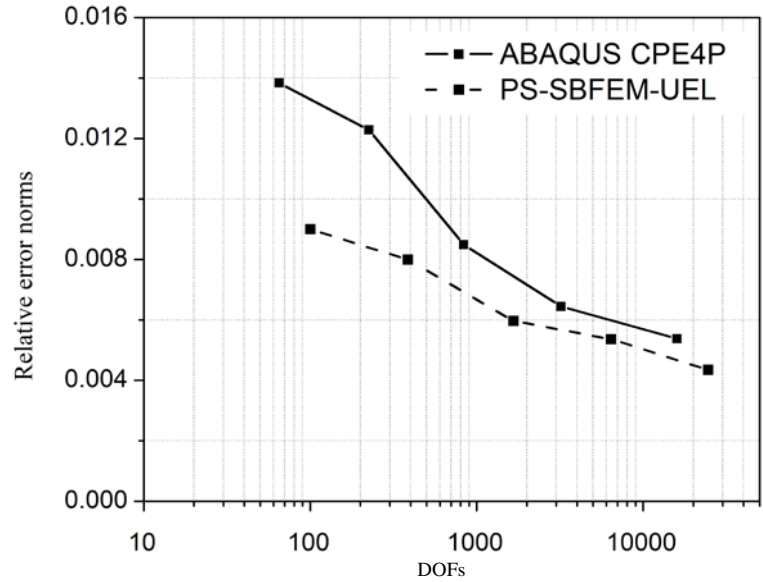
Figure 7. Comparison of convergence rate in the water head

5.2 Steady-state seepage analysis in the permeable materials

In order to show the flexibility of the PS-SBFEM using the quadtree mesh, we solve a steady-state seepage problem of permeable materials. The geometry is shown in Figure 8 (a), in which the interior of the permeable material contains impermeable material. The length and width of the permeable material are 1m, respectively. The coefficient of permeability is $k_x = k_y = 5 \times 10^{-4}$ cm/s. The quadtree mesh is shown in Figure 8 (b). The quadtree mesh has the same mesh size at the junction of the two materials, and the mesh transition area can be effectively processed without further manual intervention.

Moreover, the impermeable material does not divide the mesh, and only the impermeable boundary is set at the junction with the permeable material. To verify the accuracy of the quadtree mesh calculation, we compare the results of the Abaqus CP4EP element with similar degrees of freedom. The degrees of freedom of quadtree and CPE4P element is 11447 and 11749, respectively.

Figure 9 shows the comparison between the PS-SBFEM quadtree element and the Abaqus CPE4P element in the water head. The relative error of left edge and right edge is 0.28% and 0.32%, respectively. Furthermore, the distributions of water head obtained from the PS-SBFEM and the FEM are illustrated in Figure 10. It can be observed that the contour plots present a good agreement. Therefore, these results demonstrate PS-SBFEM accuracy and reliability for the quadtree mesh.

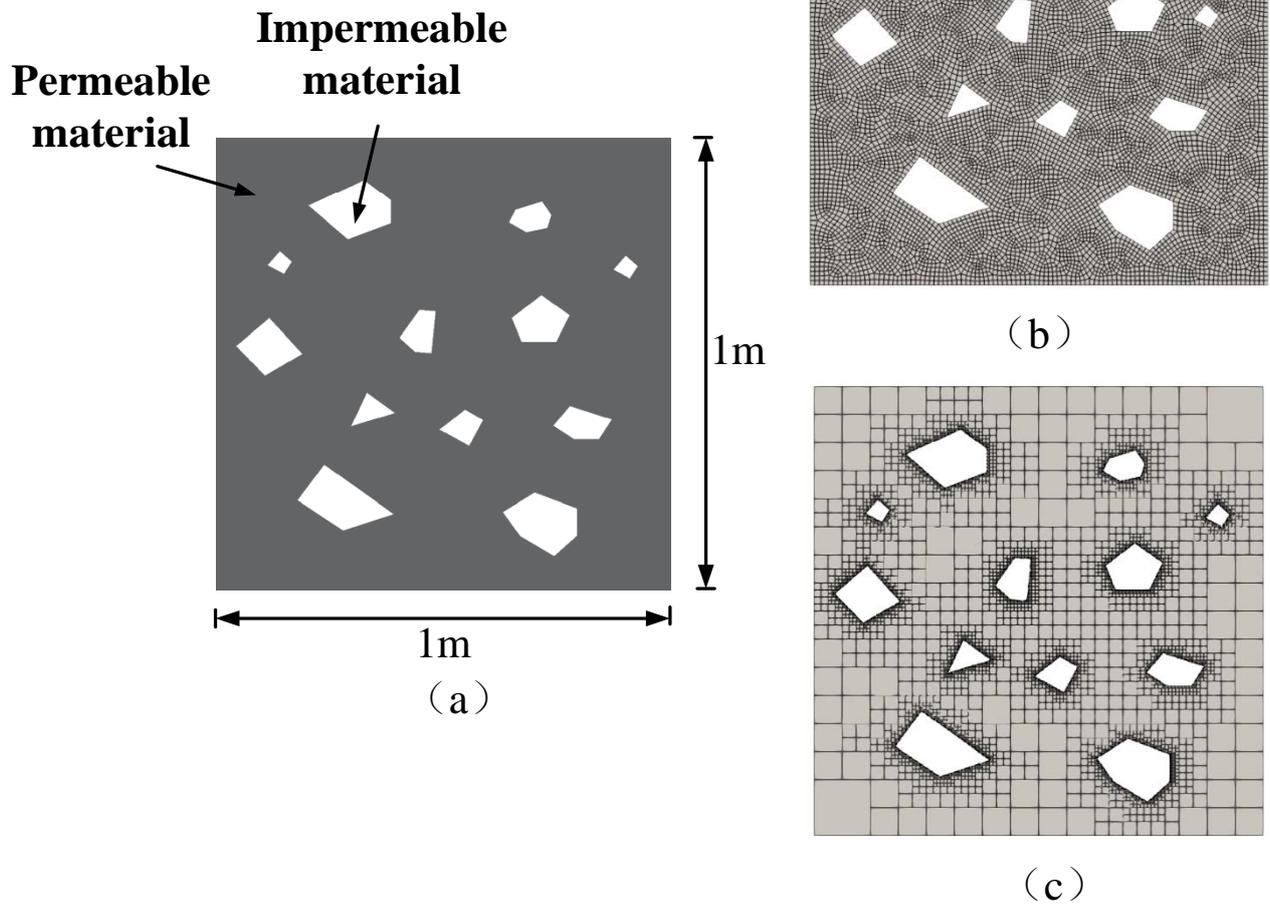

Figure 8. The permeable material geometric model and quadtree mesh; (a) geometric model; (b) Abaqus CPE4P mesh; (c) the PS-SBFEM quadtree mesh.

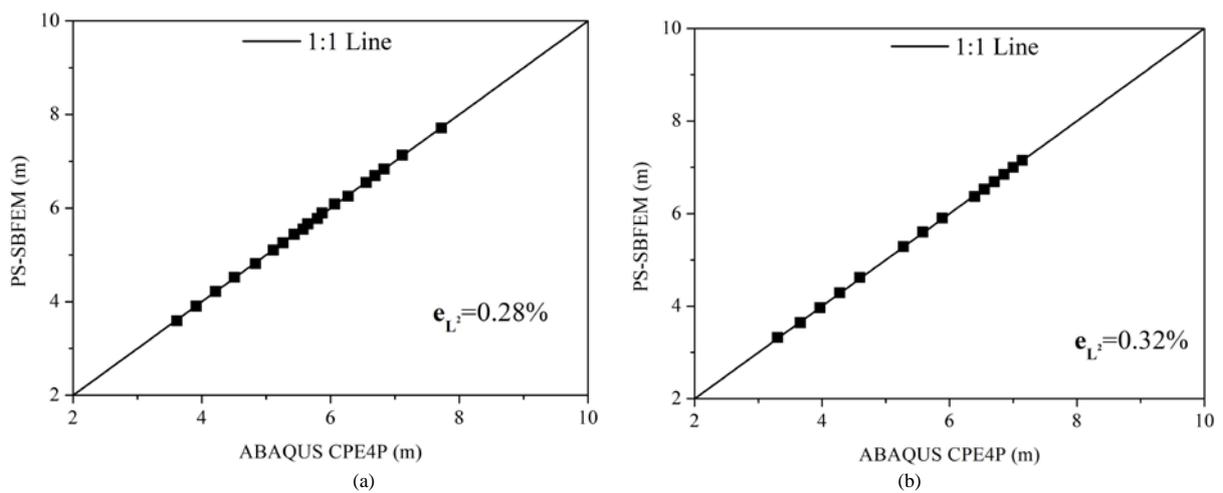

Figure 9. Comparison between the PS-SBFEM and the FEM in the water head; (a) the left edge; (b) the right edge

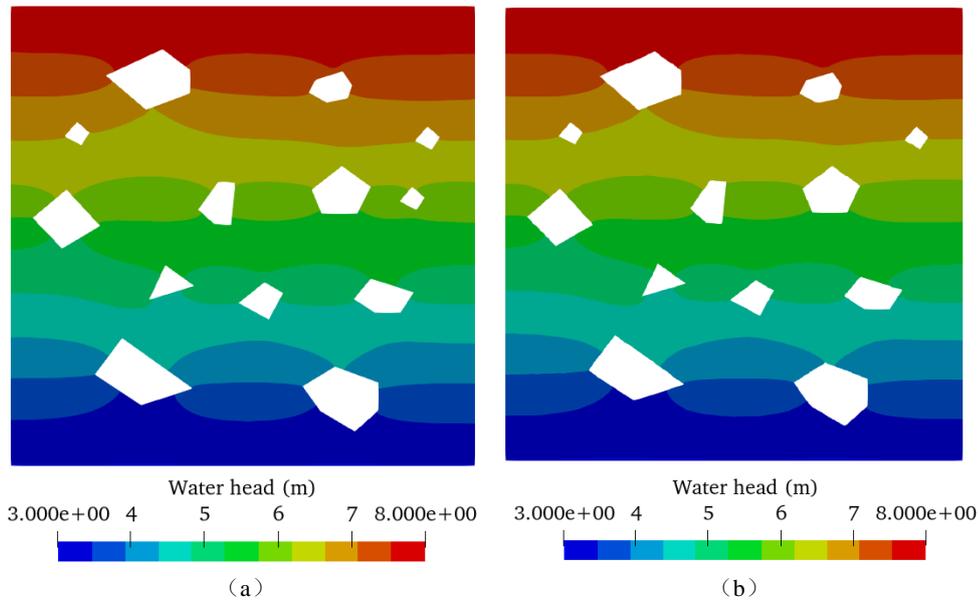

Figure 10. The water head distribution; (a) Abaqus CPE4P; (b) PS-SBFEM quadtree element.

5.3 Transient seepage analysis for the complex geometry

To demonstrate the ability of PS-SBFEM to solving complex geometry in the transient seepage, we consider a square plate ($L = 4.0m$) with a Stanford bunny cavity [38,39], as shown in Figure 11. The coefficient of permeability in $x$ and $y$ directions are considered $k_x = k_y = 5 \times 10^{-6} m/s$. The value of $S_s$ is 0.001 m$^{-1}$. As shown in Figure 11, we chose four monitor points A, B, C, D to compare results between the FEM and PS-SBFEM. The water head at the top and bottom boundaries is specified as 10m and 3m, respectively. The total time 2000 $s$, and the time step $\Delta t = 1s$ are used in the PS-SBFEM and the FEM. The PS-SBFEM and the FEM are modeled using the quadtree element and CPE4P element, respectively. As shown in Figure 12, both approaches using the same element size.

Figure 13 illustrates the history of the water head for four monitor points. The PS-SBFEM obtained solutions are in excellent agreement with the FEM for all points. When the time is greater than the 1500s, the water head of all points becomes stable. In addition, the water head of four monitor points at different times is presented in Table. 2, which shows that the relative error of four nodes is less than 1.6%. It is noted that the relative error reduces as the increment of time. Furthermore, Figure 14 shows that the distribution of water head at different times. The water head distributions are virtually the same for the FEM and PSBFEM. Therefore, the PS-SBFEM with quadtree meshes shows a good effect for solving complex geometric in the transient seepage problem.

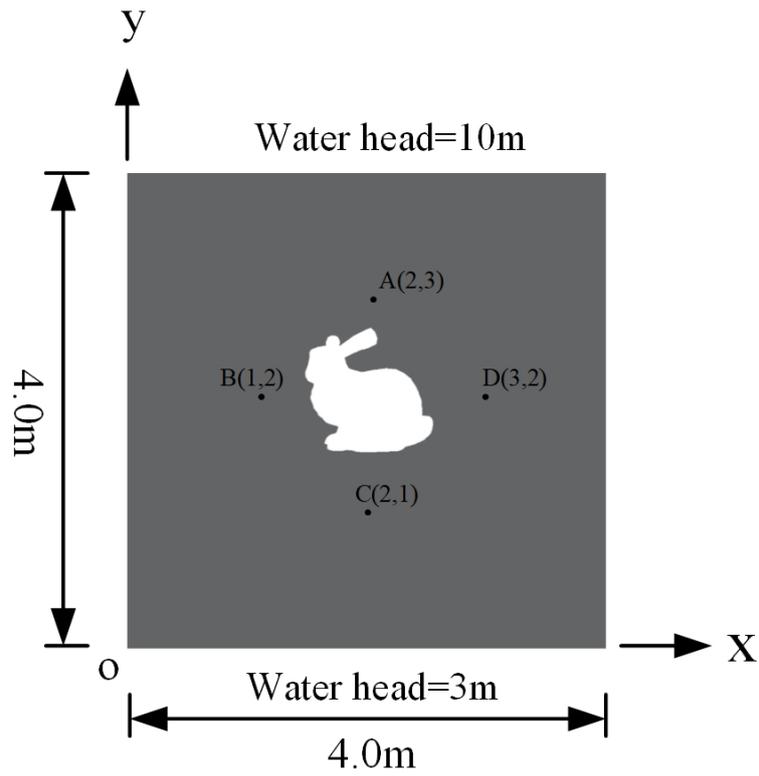

Figure 11. The geometry and boundary conditions of square plate with a Stanford bunny cavity.

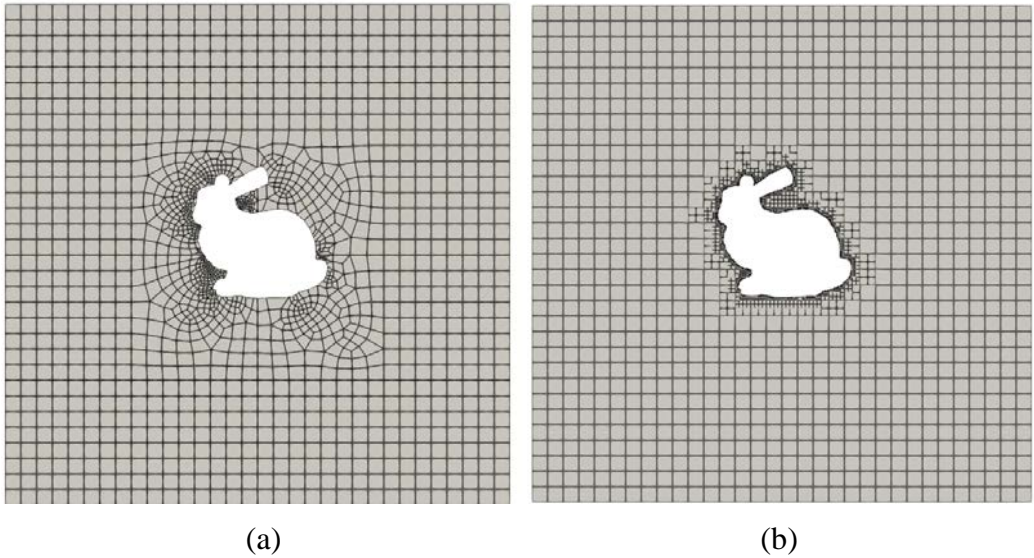

(a)                                 (b)

Figure 12. The meshes of square plate with a Stanford bunny cavity; (a) Abaqus mesh (b) quadtree mesh.

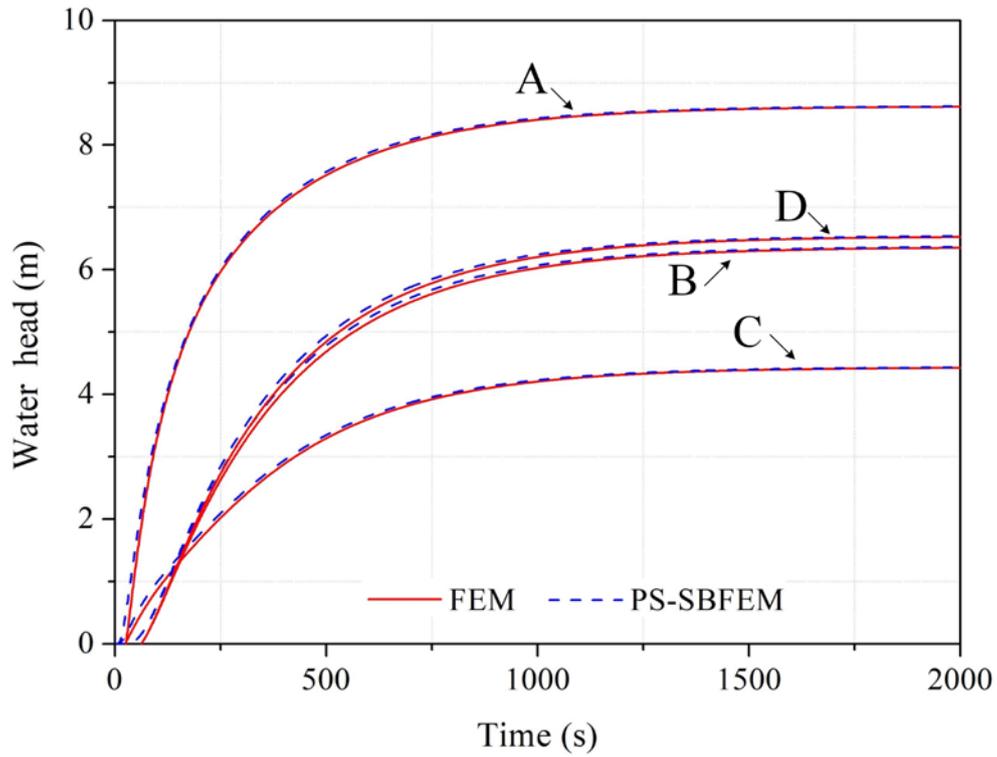

Figure 13. Comparison of the water head history of four monitor points for the FEM and PS-SBFEM

Table. 2 The water head of monitor points at the different times

| Time (s) | A FEM (m) | A PS-SBFEM (m) | B FEM (m) | B PS-SBFEM (m) | C FEM (m) | C PS-SBFEM (m) | D FEM (m) | D PS-SBFEM (m) | Relative error (%) |
|---|---|---|---|---|---|---|---|---|---|
| 500 | 7.517 | 7.571 | 4.685 | 4.783 | 3.295 | 3.348 | 4.845 | 4.940 | 1.57 |
| 1000 | 8.403 | 8.428 | 6.028 | 6.070 | 4.205 | 4.228 | 6.203 | 6.244 | 0.55 |
| 1500 | 8.578 | 8.588 | 6.297 | 6.317 | 4.388 | 4.399 | 6.472 | 6.490 | 0.24 |
| 2000 | 8.612 | 8.618 | 6.351 | 6.364 | 4.425 | 4.431 | 6.525 | 6.537 | 0.15 |

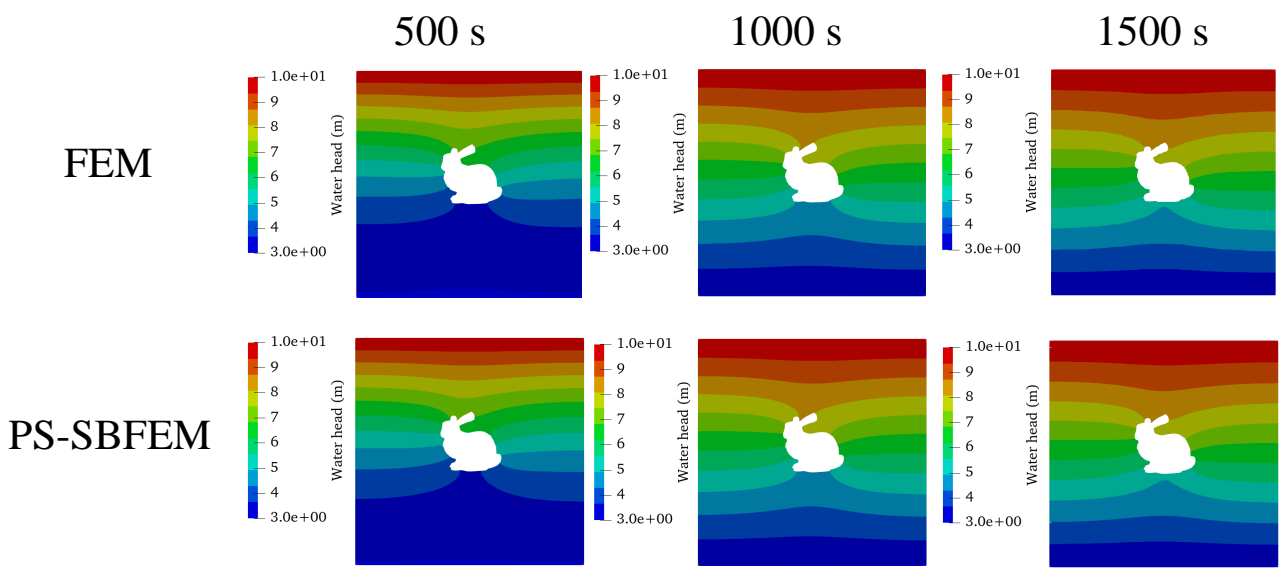

Figure 14. The water head distribution at different times using the FEM and the PS-SBFEM.

5.4 Transient seepage analysis in a concrete dam with an orthotropic foundation

In the last example, PS-SBFEM is applied to simulate the concrete dam with an orthotropic foundation. The geometry and boundary conditions are shown in Figure 15 (a). The initial water level of upstream and downstream is 10 m and 5 m, respectively. Moreover, the water level in the upstream of reservoir increased linearly to 30 m over 100 days, as illustrated in Figure 16. The material properties of $k_x$, $k_y$ and $S_s$ are 0.001 m/day, 0.0005 m/day, and 0.001 m$^{-1}$. The quadrilateral and polygonal mesh is used in the Abaqus and PS-SBFEM, respectively, as shown in Figure 15 (b) and (c). The mesh size is 5 m.

We chose a monitor point at the bottom dam to compare the PS-SBFEM and FEM results, as shown in Figure 15 (a). Table. 3 shows that the water head of monitor points at six different times. When the time is 500 days, the relative error of water head is 2.28%. However, the relative error decreases with time increasing. It can be observed that the relative error is 0.08% when the time is 3000 days. The history of water head is shown in Figure 17, where the results obtained by the two methods correspond well. It is noted that the water head becomes stable when the time is more than 2000 days. Moreover, Figure 18 illustrates the distribution of the water head at different times using the PS-SBFEM and the FEM. The two methods of results are in excellent agreement.

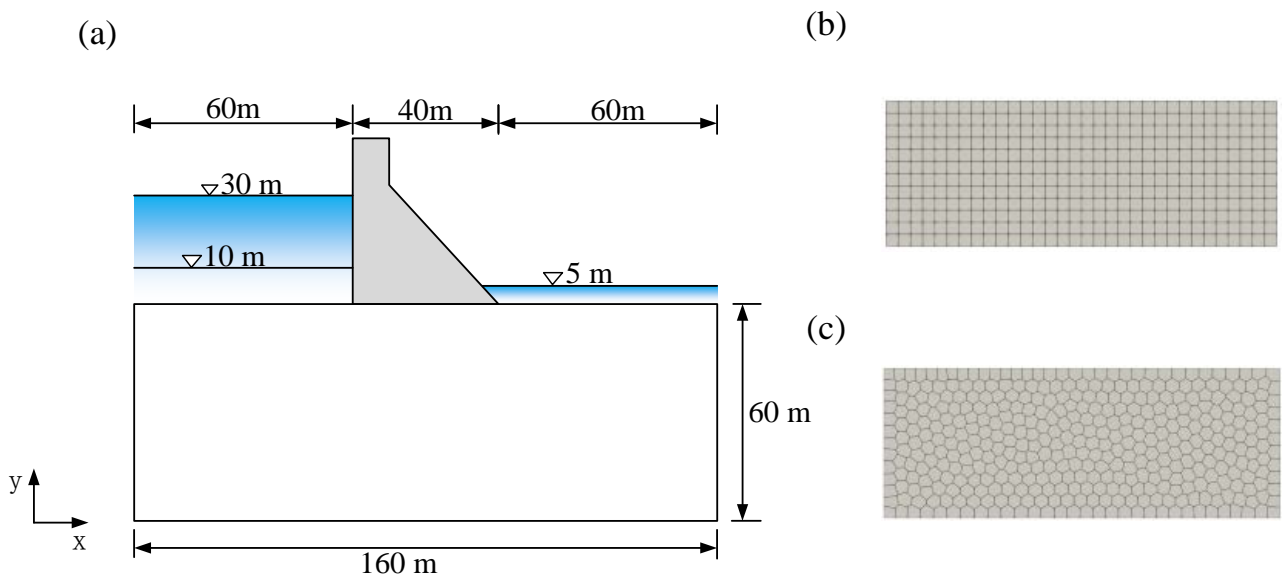

Figure 15. Transient seepage under a concrete dam constructed on an anisotropic soil; (a) the geometry and boundary conditions; (b) the FEM mesh; (c) the PS-SBFEM mesh.

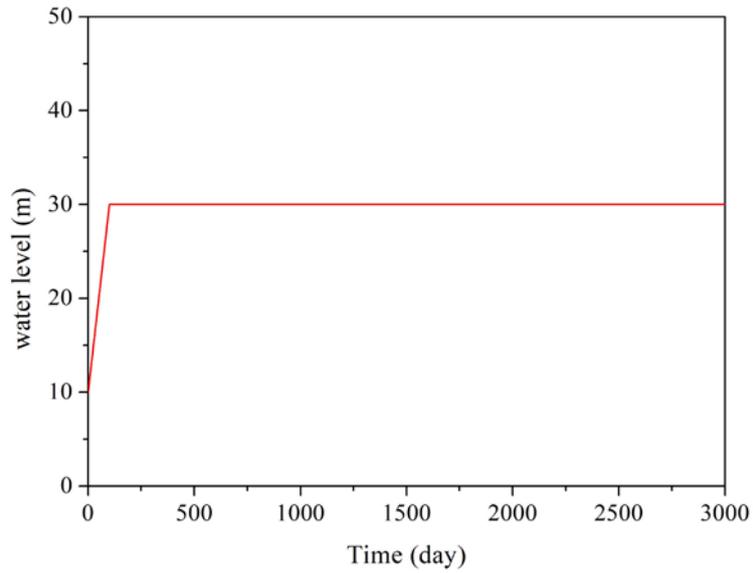

Figure 16. Variation of head with time in upstream of a concrete dam constructed on an anisotropic soil.

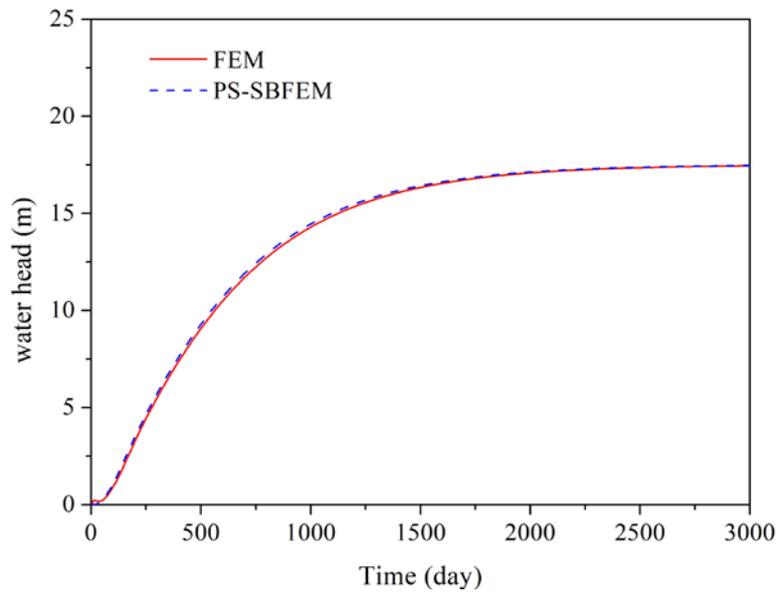

Figure 17. Comparison between the FEM and PS-SBFEM the water head history of monitor point

Table. 3 The water head of monitor points at the different times

| Time (day) | Monitor point FEM (m) | Monitor point PS-SBFEM (m) | Relative error (%) |
|---|---|---|---|
| 500 | 9.076 | 9.287 | 2.28 |
| 1000 | 14.295 | 14.454 | 1.10 |
| 1500 | 16.340 | 16.428 | 0.54 |
| 2000 | 17.087 | 17.132 | 0.26 |
| 2500 | 17.354 | 17.377 | 0.13 |
| 3000 | 17.448 | 17.462 | 0.08 |

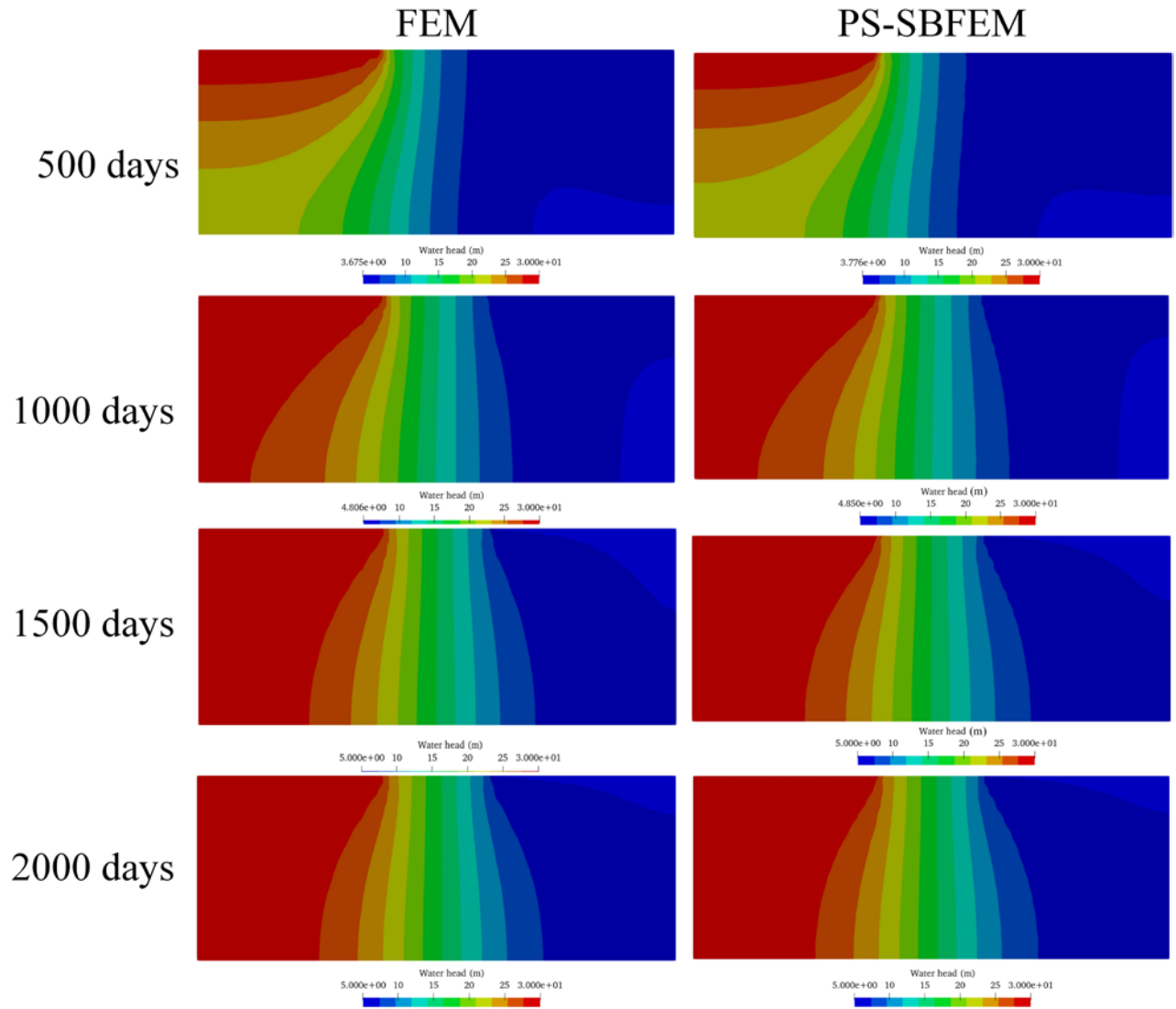

Figure 18. Comparison between the FEM and PS-SBFEM the water head distribution at different times.

## 6. Conclusion

In this work, we propose a polygonal seepage scaled boundary finite element method (PS-SBFEM) by integrating the scaled boundary finite element method (SBFEM) and the polygonal mesh technique. The implementation of PS-SBFEM is validated against the FEM by solving a several benchmark problems.

For the steady-state problems, the polygonal element of PS-SBFEM has a higher accuracy rate than the standard FEM element in the same element size. The PS-SBFEM converges to an analytical solution with an optimal convergence rate. For the transient problems, the results between PS-SBFEM and the FEM are in excellent agreement. Furthermore, the PS-SBFEM with quadtree meshes shows a good effect for solving complex geometric in the seepage problem. Hence, the proposed method is robust accurate for solving the steady-state and transient seepage problems.